\documentclass[12pt,reqno]{amsart}

\usepackage[margin=1in]{geometry}

\title[Using IFS to Reveal Biases in the Distribution of Prime Numbers]{Using Iterated Function Systems to Reveal Biases in the Distribution of Prime Numbers}

\author{Harlan J. Brothers}
\address{Brothers Technology, LLC\newline
\hphantom{00}1204 Main Street\newline
\hphantom{00}Branford, CT 06405}
\email{harlan@brotherstechnology.com}

\usepackage{placeins}
\usepackage[toc,page]{appendix}
\usepackage{bold-extra}

\usepackage[usenames]{color}
\usepackage[colorlinks=true,
linkcolor=webgreen,
filecolor=webbrown,
citecolor=webgreen]{hyperref}
\definecolor{webgreen}{rgb}{0,.5,0}
\definecolor{webbrown}{rgb}{.6,0,0}
\definecolor{ifsgreen}{rgb}{0, 0.8, 0.3}
\definecolor{ifsred}{rgb}{.9, 0.2, 0}
\definecolor{ifsblue}{rgb}{0, 0.5, 0.8}
\usepackage{gensymb}
\usepackage{paralist}
\usepackage{epsfig}
\usepackage{epstopdf}
\usepackage{amsfonts,amsmath,amssymb,amsthm}
\usepackage{comment}
\usepackage{multicol}
\usepackage{breqn}
\usepackage{comment}
\usepackage{chngpage}

\usepackage[font=footnotesize,labelfont=bf]{caption}

\begin{document}

\begin{abstract}  Iterated function systems (IFS) can be a surprisingly useful tool for studying structure in data. Here we present results stemming from a 2013 computational study by the author using IFS.  The results include fractal patterns that reveal ``repulsive" phenomena among primes in a wide range of classes, having specified arithmetic or congruence properties. Some of the phenomena shown in our computations relate to the recent, groundbreaking work of Lemke Oliver and Soundararajan on biases between consecutive primes. We do not have asymptotics to explain our results, but provide graphs, data, and detailed explanations of the phenomena.\end{abstract}

\maketitle
\vspace*{-\baselineskip}
\begin{center}
\footnotesize
(Draft, December 18, 2016)
\end{center}

\section{Introduction}

\noindent 
Graphical patterns in the primes have, for the most part, taken the form of variations of the Sieve of Eratosthenes or of Ulam's spiral \cite{Weisstein}.  In this paper, we present results that were initially generated during a 2013 computational study using \emph{Mathematica}.  The study examined relations between various arithmetic classes of prime numbers using ideas from iterated function systems (IFS).  IFS provides a method for the generation and analysis of fractals, a field of deep interest to the author \cite{Brothers1}.

\vskip 6pt

IFS was first popularized by Michael Barnsley, a pioneer in the field \cite{Barnsley1}.  In his words, iterated function systems ``provide models for certain plants, leaves, and ferns, by virtue of the self-similarity which often occurs in branching structures in nature'' \cite{Barnsley2}.
Similar ``branching structures,'' such as divisor lattices, arise naturally in classical number theory. The author wondered if IFS could be used to uncover such patterns in the distribution of primes.  The results were intriguing and unexpected.  The author noticed the appearance of fractal patterns in these IFS graphs generated from prime numbers.  The patterns revealed a puzzling ``repulsion'' phenomena between consecutive primes across congruence classes that indicate a non-uniform distribution.

\vskip 6pt

The author shared these early results with Michael Frame at Yale University and Michael Barnsley, neither of whom could explain the unexpected repulsive phenomena \cite{Frame-Brothers, Barnsley-Brothers}. He then put the project aside due to the inherent difficulty in establishing the mechanism for these visual relationships. Indeed, classical facts like Fermat's little theorem, Dirichlet's theorem on primes in arithmetic progressions, as well as modern realizations by Zhang, Maynard, Tao, et al. related to small gaps between primes, are all, in some sense, woven together graphically in images that follow.

\vskip 6pt

Recently, the existence of patterns like these have been explained in some detail by Robert Lemke Oliver and Kannan Soundararajan in their groundbreaking analysis \cite{LO-Sound}.

\vskip 8pt

\section{Driven IFS} 
\noindent
Broadly speaking, a fractal generated by IFS is defined by a collection of transformations.  These transformations take the form of contraction maps over a complete metric space and can be grouped into two broad categories: deterministic IFS and random IFS.

\vskip 6pt

For deterministic IFS and contraction maps $\{\text{T}_{i}:\mathbb R^{2}\to\mathbb R^{2}\mid i=1,2,\dots ,N\},N \in \mathbb{N}$, the collage map \emph{T} is defined on the set $K(\mathbb R^{2})$ of compact subsets of $\mathbb R^{2}$ by
 
\begin{equation}
T(C) = \textmd{T}_{1}(C) \cup \dots  \cup \textmd{T}_{N}(C)
\label{eq1}
\end{equation}
\vskip 6pt

\noindent where $\text{T}_{i}(C) = \{\text{T}_{i}(x, y): (x, y) \in C\}$.  There is then a unique compact set $A$ satisfying $T(A) = A$ and for any compact set $B$,

\begin{equation}
\lim_{k \to \infty}T^k({B})={A}
\label{eq2}
\end{equation}

\vskip 6pt

\noindent where $A$ is called the \emph{attractor} of the IFS.\footnote{For details on the derivation of Eq. (\ref{eq2}), see \cite{Frame1}.}

The set of affine transformations for scaling by $r$ in the $x$-direction and $s$ in the $y$-direction, with rotations by $\theta$ and $\phi$, and translations by $e$ and $f$ are given in matrix form by:

\[
\begin{bmatrix} 
x\\
y  
\end{bmatrix}
\rightarrow
\begin{bmatrix} 
r\operatorname{cos}(\theta) &-s\operatorname{sin}(\phi) \\
r\operatorname{sin}(\theta) &\  s\operatorname{cos}(\phi)
\end{bmatrix}
\begin{bmatrix} 
x\\
y  
\end{bmatrix}
+
\begin{bmatrix} 
{e}\\
{f}
\end{bmatrix}
\]

\vskip 6pt

\begin{figure}[htb]
\centering
\epsfxsize = 160pt
\epsfbox{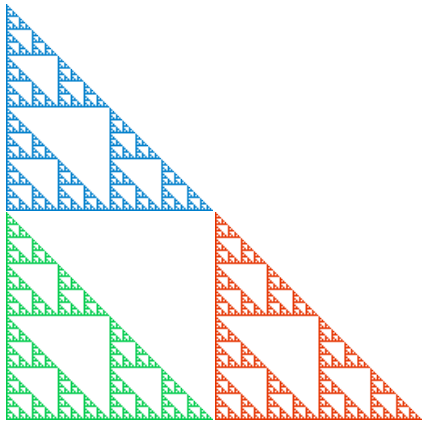}
\caption{The Sierpinski gasket decomposed into three copies of itself.}
\label{fig1}
\end{figure}

The canonical example of a deterministic IFS is the Sierpinski gasket (see Figure \ref{fig1}) for which the attractor is defined by the transformations

\begin{equation}
\begin{aligned}
\textmd{T}_{1}(x,y) &= (x/2, y/2)\\
\textmd{T}_{2}(x,y) &= (x/2, y/2) + (1/2, 0)\\
\textmd{T}_{3}(x,y) &= (x/2, y/2) + (0, 1/2).
\end{aligned}
\label{eq3}
\end{equation}

\vskip -6 pt

\noindent Table \ref{table1} shows the matrix parameters for the gasket.

\begin{table}[h!]
\footnotesize
\begin{center}
\begin{tabular}{| c || c | c | c | c | c | c |}
\hline
		& $\boldsymbol{r}$ & $\boldsymbol{s}$ & $\boldsymbol{\theta}$ & $\boldsymbol{\phi}$ & $\boldsymbol{e}$ & $\boldsymbol{f}$\\ \hline
\hline
T$_{1}$ & \color{ifsgreen} {.5} & \color{ifsgreen}{.5} & \color{ifsgreen}{0} & \color{ifsgreen}{0} & \color{ifsgreen}{0} & \color{ifsgreen}{0}\\ \hline
T$_{2}$ & \color{ifsred}{.5} & \color{ifsred}{.5} & \color{ifsred}{0} & \color{ifsred}{0} & \color{ifsred}{.5} & \color{ifsred}{0}\\ \hline
T$_{3}$ & \color{ifsblue}{.5} & \color{ifsblue}{.5} & \color{ifsblue}{0} & \color{ifsblue}{0} & \color{ifsblue}{0} & \color{ifsblue}{.5}\\ \hline
\end{tabular}
\caption{Matrix parameters for the three transformations describing the gasket, color coded to Figure \ref{fig1}.}
\label{table1}
\end{center}
\end{table}

\vskip 4 pt

Random IFS uses the same transformations as deterministic IFS, however, rather than applying all T$_{i}$ simultaneously and then iterating on the output, we assign a probability $p_{i}$ to each T$_{i}$ where $\sum_{i=1}^{N}p_i=1$.  In the simplest case, the $p_{i}$ are equal -- all T$_{i}$ are applied with equal likelihood. This is equivalent to deterministic IFS.  We can use the Sierpinski gasket to illustrate this method.

\vskip 6pt

First, choose  a fixed point $(x_{0}, y_{0})$ from one of the T$_{i}$ in (\ref{eq3}).  Using a uniformly distributed random sequence of $k$ numbers $\{n_{1}, n_{2}, \ldots\, n_{k}\}$, $n \in \{1, \ldots, N\}$, we generate a sequence of points:

\begin{equation}
\begin{aligned}
&(x_{1}, y_{1}) = \textmd{T}_{n_{1}}(x_{0}, y_{0})\\
&(x_{2}, y_{2}) = \textmd{T}_{n_{2}}(x_{1}, y_{1})\\
& \hskip 16pt \vdots \hskip 52pt \vdots\\
&(x_{k}, y_{k}) = \textmd{T}_{n_{k}}(x_{k-1}, y_{k-1})
\end{aligned}
\label{eq4}
\end{equation}

\noindent This sequence of points will eventually fill in the gasket to any desired resolution.  For further details regarding theory and variations of random IFS, see \cite{Barnsley3, Barnsley4, Barnsley5}.

\vskip 6pt

To explore patterns in the distribution of primes, we will use a type of random IFS called \emph{driven} IFS \cite{Frame2}.  To start, consider what happens if we add a fourth transformation to the set of gasket transformations in (\ref{eq3}) that shrinks everything by 1/2 and translates it up and to the right 1/2 unit:

\begin{table}[htbp]
\footnotesize
\begin{center}
\begin{tabular}{| c || c | c | c | c | c | c |}
\hline
		&$\boldsymbol{r}$ & $\boldsymbol{s}$ & $\boldsymbol{\theta}$ & $\boldsymbol{\phi}$ & $\boldsymbol{e}$ & $\boldsymbol{f}$\\ \hline
\hline
T$_{1}$ & .5 & .5 & 0 & 0 & 0 & 0\\ \hline
T$_{2}$ & .5 & .5 & 0 & 0 & .5 & 0\\ \hline
T$_{3}$ & .5 & .5 & 0 & 0 & 0 & .5\\ \hline
T$_{4}$ & .5 & .5 & 0 & 0 & .5 & .5\\ \hline
\end{tabular}
\caption{Matrix parameters for the filled-in unit square.}
\label{table2}
\end{center}
\end{table}

\vskip -10pt

\begin{figure}[htbp]
\centering
\epsfxsize = 160pt
\epsfbox{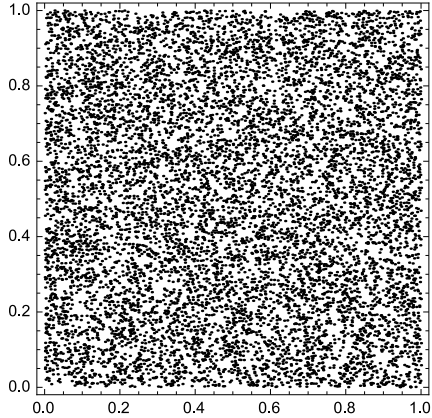}
\caption{Plot of 1000 points generated by randomly applying the transforms in Table (\ref{table2}).}
\label{fig2}
\end{figure}

\noindent If we apply these transforms randomly (and uniformly) to a starting point $(x_{0}, y_{0})$, we generate the filled-in unit square (see Figure \ref{fig2}).

\vskip 6pt

One might reasonably ask what happens if the T$_{i}$ in Table \ref{table2} are not applied randomly? Indeed, ``driving" this IFS with data can reveal unexpected patterns to the extent that the data deviates from a random sequence.  Perhaps the earliest example of this approach was the work of H. Joel Jeffrey in exploring representations of gene structure \cite{Jeffrey}.

\vskip 6pt

Until recently, it was assumed that the residues for a given modulus $q$ are evenly distributed across the primes.  As Section 4 demonstrates, driven IFS reveals that this not the case.

\vskip 16pt

\section{IFS Addresses} 
\begin{figure}[htb]
\centering
\epsfxsize = 80pt
\epsfbox{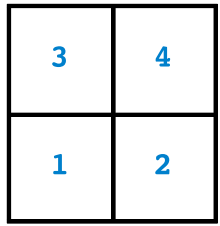}
\caption{IFS length 1 addresses.}
\label{fig3}
\end{figure}

\noindent 
To understand the dynamics of plots associated with driven IFS, it is important to understand the standard address system for IFS on the unit square $S$. The transformations T$_{1}$, T$_{2}$, T$_{3}$, and T$_{4}$ in Table \ref{table2} can be thought of as moving a point halfway toward the vertices at $(0, 0), (1, 0), (0, 1)$, and $(1, 1)$, respectively.  We associate  the quadrants for each vertex with a length 1 address as shown in Figure \ref{fig3}, so that $S$ is decomposed into:

\vskip 8pt

\begin{equation}
{S} = \textmd{T}_{1}({S}) \cup \textmd{T}_{2}({S}) \cup \textmd{T}_{3}( {S})  \cup \textmd{T}_{4}({S})
\label{eq5}
\end{equation}

\vskip 6pt

\vskip 6pt

\begin{figure}[htb]
\centering
\epsfxsize = 120 pt
\epsfbox{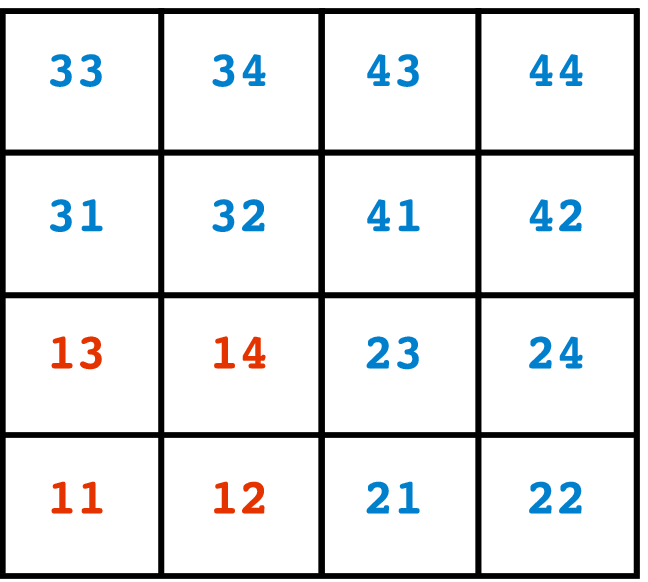}
\caption{IFS length 2 addresses. T$_{1}$ addresses are red.}
\label{fig4}
\end{figure}

The quadrant represented by each length 1 address $i$ can in turn be divided by iterating the decomposition process to obtain a length 2 address $ij$. Looking at quadrant 1 in Figure \ref{fig4}, we have:

\begin{equation}
\textmd{T}_{1}(S) = \textmd{T}_{1}\textmd{T}_{1}(S) \cup \textmd{T}_{1}\textmd{T}_{2}(S) \cup \textmd{T}_{1}\textmd{T}_{3}(S) \cup \textmd{T}_{1}\textmd{T}_{4}(S)
\label{eq6}
\end{equation}

\vskip 8pt

\noindent The other quadrants are similarly subdivided. This iterative decomposition continues to whatever resolution is required to specify a location.

\vskip 6pt

It is important to note that addresses are read from right to left.  Just as with the order of composition of functions, $ij$ is the address of T$_{i}$T$_{j}(S)$.  That is, the left-most digit is the index of the most recently applied transformation.\footnote{For more information on addresses in IFS, see \cite{Frame3}.}

\vskip 16pt

\section{Repulsive Behavior in the Primes} 
\noindent 
Here we will use the notation $\pi(x_{0},x; q, (a, b))$ to denote the number of consecutive primes beginning with $x_{0}$ up to $x$ that are members of the reduced residue classes $a$ (mod $q$) and $b$ (mod $q$). 

\vskip 6pt

We wish to look for deviations from the uniform distribution shown in Figure \ref{fig2}. It makes sense, then, to examine the primes using primitive residue classes modulo $n$ for $n \in \{5,8,10,12\}$, all of which have 4 elements.  For our purposes, the groups $M_{5}=\{1,2,3,4\}$ and $M_{10}=\{1,3,7,9\}$ are equivalent given that the elements of $M_{5}$ generate the same four sequences of residues over the primes as do those of $M_{10}$ (i.e., $7\equiv 2~(\mathrm{mod}~5)$ and $9\equiv 4~(\mathrm{mod}~5)$).  We use $M_{10}$ for the sake of convenience.

\vskip 6pt

We therefore consider $M_{8}$, $M_{10}$, and $M_{12}$.  For each of these, there are $4!$ ways to assign residues to the vertices of our unit square.  Out of the 24 possible arrangements for each group, we can eliminate the permutations that are equivalent under rotation and reflection (these will simply rotate or flip our IFS plot).  Table \ref{table1} shows the remaining three inequivalent permutations for each of the three groups.

\begin{table}[h!]
\footnotesize
\begin{center}
\begin{tabular}{| c | c | c |}
\hline
$\boldsymbol{M_{8}}$ & $\boldsymbol{M_{10}}$ & $\boldsymbol{M_{12}}$\\ \hline
\hline
\{1, 3, 5, 7\} & \{1, 3, 7, 9\} & \{1, 5, 7, 11\} \\ \hline
\{1, 3, 7, 5\} & \{1, 3, 9, 7\} & \{1, 5, 11, 7\} \\ \hline
\{1, 5, 3, 7\} & \{1, 7, 3, 9\} & \{1, 7, 5, 11\}\\ \hline
\end{tabular}
\caption{Free circular permutations of order 4 for modulo multiplication groups $M_{m}$.}
\label{table3}
\end{center}
\end{table}

\vskip -12pt

Of the free circular permutations listed in Table \ref{table3}, we can eliminate the final row from our examination given that the corresponding IFS is simply a 90$^{\circ}$ counterclockwise rotation, followed by a horizontal reflection of the IFS generated by row 1 (see Figure \ref{fig3}).

\vskip 6pt

\begin{figure}[htbp]
\centering
\epsfxsize = 420pt
\epsfbox{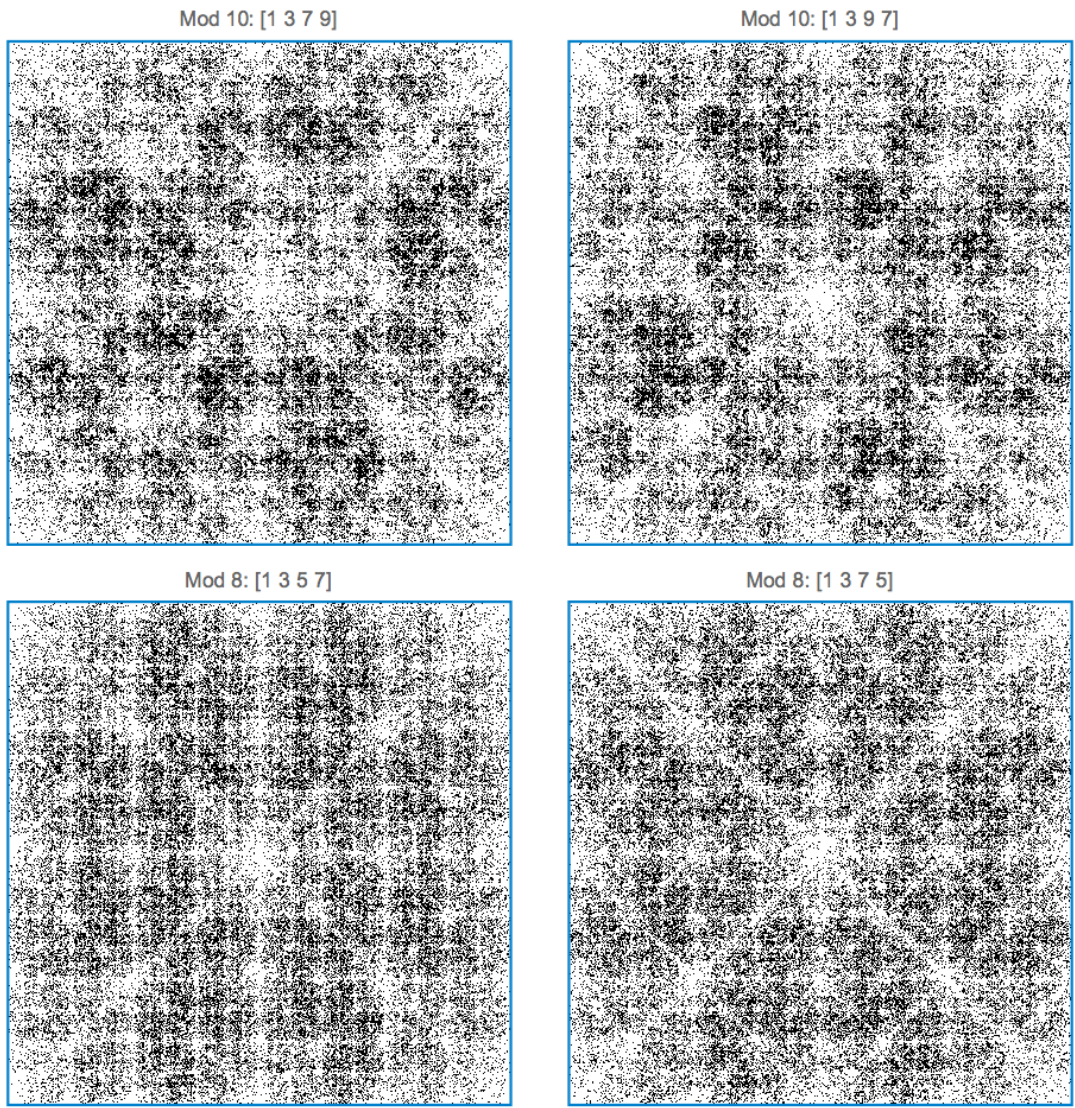}
\caption{IFS plots for the residues modulo 10 and modulo 8 for primes $p$ up to $10^{6}$, beginning with $p=7$.}
\label{fig5}
\end{figure}

For any four values $a, b, c, d$, we denote the mapping of $a\rightarrow$ T$_{1}$, $b\rightarrow$ T$_{2}$, $c\rightarrow$ T$_{3}$, and $d\rightarrow$ T$_{4}$, by $\left[a~b~c~d\right]$.  The top of Figure \ref{fig5} shows $10^{6}$ primes $p$ modulo 10, beginning with $p=7$.  On the left side is $\left[1~3~7~9\right]$ and on the right is $\left[1~3~9~7\right]$.\footnote{Higher resolution versions of these and other plots are available at: (\emph{needs link})\\}

\vskip 6pt

\begin{figure}[htbp]
\centering
\epsfxsize = 280pt
\epsfbox{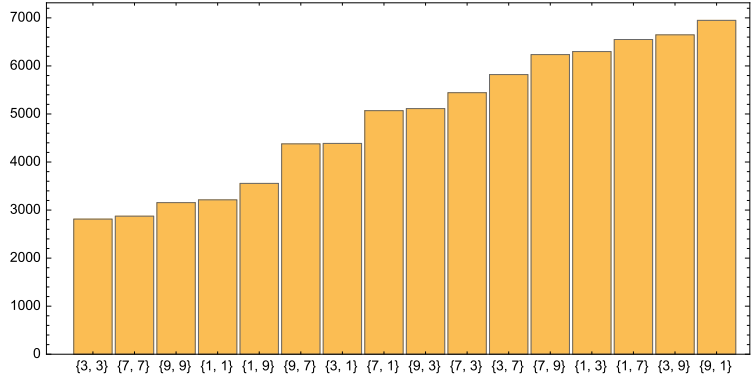}
\caption{Histogram of residues modulo 10 corresponding to length 2 addresses for primes $p$ up to $10^{6}$, beginning with $p=7$. See Figure \ref{fig5}.}
\label{fig6}
\end{figure}

In this case, the two plots are very similar.  Their defining feature is that the corners of the unit square are relatively sparsely populated.  This tells us that the pair of transformations T$_{i}$T$_{j}$ occurs less frequently for $i=j$.  The histogram in Figure \ref{fig6} shows the frequencies of length 2 addresses and indeed corroborates this observation.  Within this range, the least frequent pairs occur less than half as often as the most frequent.

\begin{figure}[htbp]
\centering
\epsfxsize = 280pt
\epsfbox{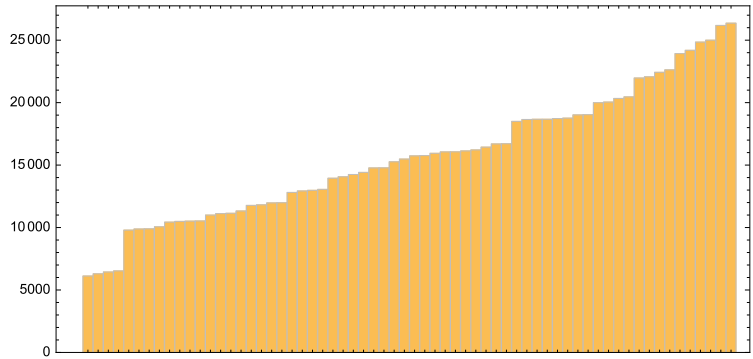}
\caption{Histogram of residues modulo 10 corresponding to length 3 addresses for primes $p$ up to $10^{6}$, beginning with $p=7$. The sequence \{7, 7, 7\} appears least frequently and the sequence \{9, 1, 7\} appears most frequently.  See Figure \ref{fig5}.}
\label{fig7}
\end{figure}

\vskip 6pt

By extension, T$_{i}$T$_{j}$T$_{k}$ should occur even less often as a percentage of length 3 addresses for $i=j=k$.  The histogram in Figure \ref{fig7} shows that within this range, the least frequent triples occur roughly one quarter as often as the most frequent. Indeed, the sparsely populated corners make a recursive appearance at each scale as defined by the address length.  This gives a grid-like quality to the plot (and facilitates identifying address locations without actual gridlines).

\vskip 6pt

Specifically, this means that a prime ending in 1 is not as likely to be immediately followed by another prime ending in 1.  The same can be said for those ending in 3, 7, or 9.  The residue classes tend to avoid appearing adjacently.

\vskip 6pt

Visual inspection can tell us more.  Looking at the modulo 10 plots at the top of Figure \ref{fig5}, we see that  address 14 on the left side, and address 13 on the right side, are darker, indicating that the region is visited more often than neighboring addresses.  This corresponds to T$_{1}$T$_{4}(S)$ on the left and T$_{1}$T$_{3}(S)$ on the right.  Recalling the order of operations, this means that a prime ending in 1 is most likely to follow a prime ending in 9.  The histogram in Figure \ref{fig6} confirms this observation.

\vskip 6pt

The IFS plots for $M_{8}$ and $M_{12}$ show similar results, though with somewhat different patterns.  For comparison, Figure \ref{fig5} includes the IFS plots for $M_{8}$.  The more heterogeneous appearance of $M_{8}$ arises in part from the fact that the residues in general are more evenly distributed as seen in Table \ref{table4} ($M_{12}$ is included for reference).

\begin{table}[htbp]
\footnotesize
\begin{center}
\begin{tabular}{| c | c | c || c | c | c || c | c | c |}\hline
$\boldsymbol{M}_{8}$ $(a,b)$ & Count & ${\sigma}$ (counts) & $\boldsymbol{M}_{10}$ $(a,b)$ & Count & ${\sigma}$ (counts) & $\boldsymbol{M}_{12}$ $(a,b)$ & Count & ${\sigma}$ (counts)\\ \hline \hline
(5, 5) & 3,047 & 1091.14 & (3, 3) & 2,812 & 1452.11 & (7, 7) & 2,921 & 1412.04\\ \hline
(7, 7) & 3,107 &  & (7, 7) & 2,873 &  & (11, 11) & 2,941 & \\ \hline
(3, 3) & 3,117 &  & (9, 9) & 3,155 &  & (1, 1) & 2,944 & \\ \hline
(1, 1) & 3,128 &  & (1, 1) & 3,213 &  & (5, 5) & 2,961 & \\ \hline
(1, 5) & 5,190 &  & (1, 9) & 3,555 &  & (5, 1) & 4,365 & \\ \hline
(5, 1) & 5,239 &  & (9, 7) & 4,378 &  & (11, 7) & 4,428 & \\ \hline
(7, 3) & 5,250 &  & (3, 1) & 4,387 &  & (7, 5) & 5,018 & \\ \hline
(3, 7) & 5,276 &  & (7, 1) & 5,069 &  & (1, 11) & 5,026 & \\ \hline
(3, 1) & 5,565 &  & (9, 3) & 5,112 &  & (5, 11) & 5,227 & \\ \hline
(1, 3) & 5,614 &  & (7, 3) & 5,443 &  & (1, 7) & 5,261 & \\ \hline
(1, 7) & 5,621 &  & (3, 7) & 5,820 &  & (7, 1) & 5,269 & \\ \hline
(7, 1) & 5,621 &  & (7, 9) & 6,236 &  & (11, 5) & 5,298 & \\ \hline
(5, 7) & 5,664 &  & (1, 3) & 6,299 &  & (1, 5) & 6,334 & \\ \hline
(5, 3) & 5,672 &  & (1, 7) & 6,550 &  & (7, 11) & 6,460 & \\ \hline
(7, 5) & 5,691 &  & (3, 9) & 6,647 &  & (11, 1) & 6,987 & \\ \hline
(3, 5) & 5,695 &  & (9, 1) & 6,948 &  & (5, 7) & 7,057 & \\ \hline
\end{tabular}
\caption{Frequencies for appearances of pairs of successive residues for modulo multiplication groups $M_{q}$, sorted by count, where Count=$\pi(7,10^{6}; q, (a, b))$ and standard deviation $\sigma$ is taken across counts for $M_{q}$. See Figure \ref{fig5}.}
\label{table4}
\end{center}
\end{table}

\FloatBarrier

\section{Looking Deeper} 
\begin{figure}[htbp]
\centering
\epsfxsize = 204pt
\epsfbox{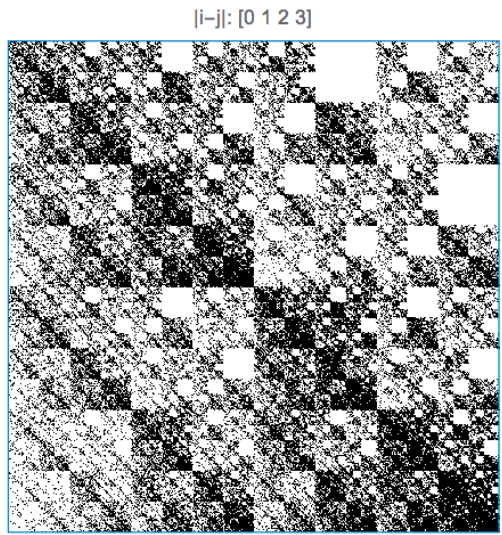}
\caption{IFS plot for the absolute differences in the transformation indices corresponding to residues modulo 10 of 200,000 primes $p$, beginning with $p=7$.}
\label{fig8}
\end{figure}

\noindent 
We can also examine the dynamics of the way in which the sequence of transformations itself changes.  One way to do this would be to look at the absolute difference of consecutive values of index $\emph{i}\in \{1, 2, 3, 4\}$ for T$_{i}$.  However, this would impose a dominant artifact in the IFS plot, unrelated to any inherent pattern, as seen in Figure \ref{fig8}.   The empty blocks occurring at 424 and 434 are referred to as ``forbidden addresses."  Indeed, here they are forbidden because the absolute differences take the values 0, 1, 2, or 3 and the only jump of length 3 must begin with either a 1 or a 4. Any absolute difference address starting with 424 or 434 simply cannot happen.

\vskip 6pt

A better approach then is to look at the forward rotational distance for $\mathbb Z_{4}$ whose distance measures, shown in Table \ref{table5}, are equally distributed.  For T$_{i}$T$_{j}$ in the previous section, keeping in mind that T$_{i}$ is applied after T$_{j}$, we denote this distance by $\|ij\|_{\textmd{rot}}~=(i-j) \pmod{4}$.

\vskip -8pt

\begin{table}[h!]
\footnotesize
\begin{center}
\begin{tabular}{| c | c | c | c | c |}
\hline
           & \textbf{1} & \textbf{2} & \textbf{3} & \textbf{4}\\ \hline
\hline
\textbf{1} & 0 & 3 & 2 & 1 \\ \hline
\textbf{2} & 1 & 0 & 3 & 2\\ \hline
\textbf{3} & 2 & 1 & 0 & 3\\ \hline
\textbf{4} & 3 & 2 & 1 & 0\\ \hline
\end{tabular}
\caption{Forward rotational distance measure $\|ij\|_{\textmd{rot}}$ for $\mathbb Z_{4}$, for row $i$ and column $j$.}
\label{table5}
\end{center}
\end{table}

\vskip -24pt

\begin{figure}[h!]
\centering
\epsfxsize = 420pt
\epsfbox{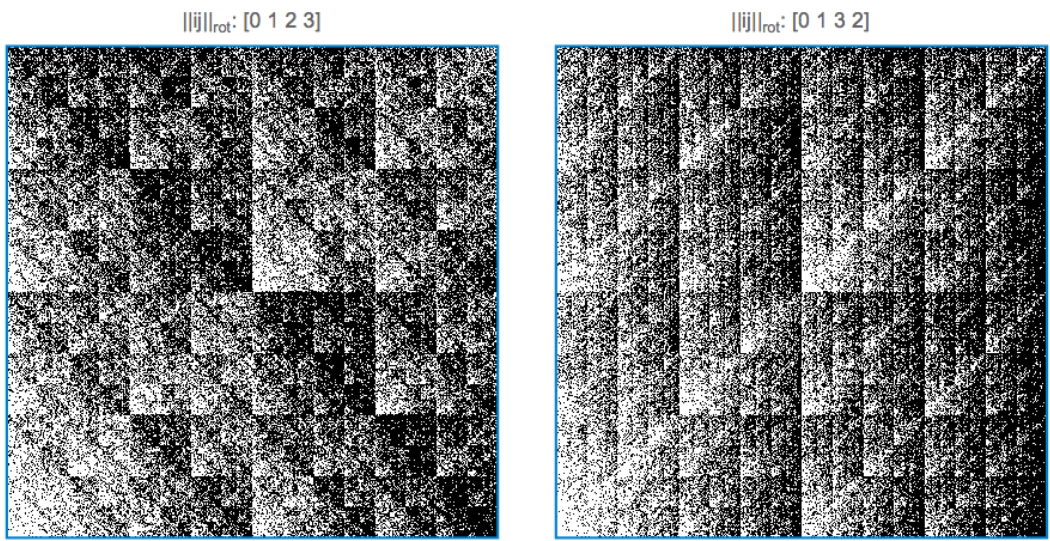}
\caption{IFS plot for the forward rotational distance $\|ij\|_{\textmd{rot}}$ between the transformation indices for the residues modulo 10 of 200,000 primes $p$, beginning with $p=7$.  Note the similarity of the left-hand plot to Figure \ref{fig8}.}
\label{fig9}
\end{figure}

\vskip 6pt

In left-hand plot of Figure \ref{fig9}, for $\|ij\|_{\textmd{rot}}~=n$, we apply T$_{n+1}$.  We can see in the lower-left corner of both plots, the light areas indicate that repeated residues, $\|ij\|_{\textmd{rot}}~=0$, occur less frequently. In the lower-right corners, the dark areas indicate that circularly increasing adjacent residues, $\|ij\|_{\textmd{rot}}~=1$, tend to occur more frequently. Furthermore, the dark diagonal band on the left from T$_{2}$ to T$_{3}$ (and corresponding vertical dark band on the border of the right) indicates $\|ij\|_{\textmd{rot}}~=1$ and $\|ij\|_{\textmd{rot}}~=2$ often follow each other. Looking at the first 10 million primes, Table \ref{table6} bears out this observation. From a dynamic standpoint, not only do the same residue classes avoid appearing in sequence, but it seems that a single ``step'' forward through the residue classes is the most likely motion between adjacent primes.  

\begin{table}[htbp]
\footnotesize
\begin{center}
\begin{tabular}{| c | c | c | c | c |}
\hline
$\mathbf{\|}$\emph{ij}$\mathbf{\|}_{\text{rot}}$ &  {Count} & {Percentage} \\ \hline
\hline
0 &  1,737,431 & 17.374\\ \hline
1 &  3,048,086 & 30.481\\ \hline
2 &  2,819,299 & 28.193\\ \hline
3 &  2,395,183 & 23.952\\ \hline
\end{tabular}
\caption{Frequencies for $\|ij\|_{\textmd{rot}}$ for $10^{7}$ primes, beginning with $p=7$.}
\label{table6}
\end{center}
\end{table}

\FloatBarrier

\section{Twin Primes}
\begin{figure}[htbp]
\centering
\epsfxsize = 420pt
\epsfbox{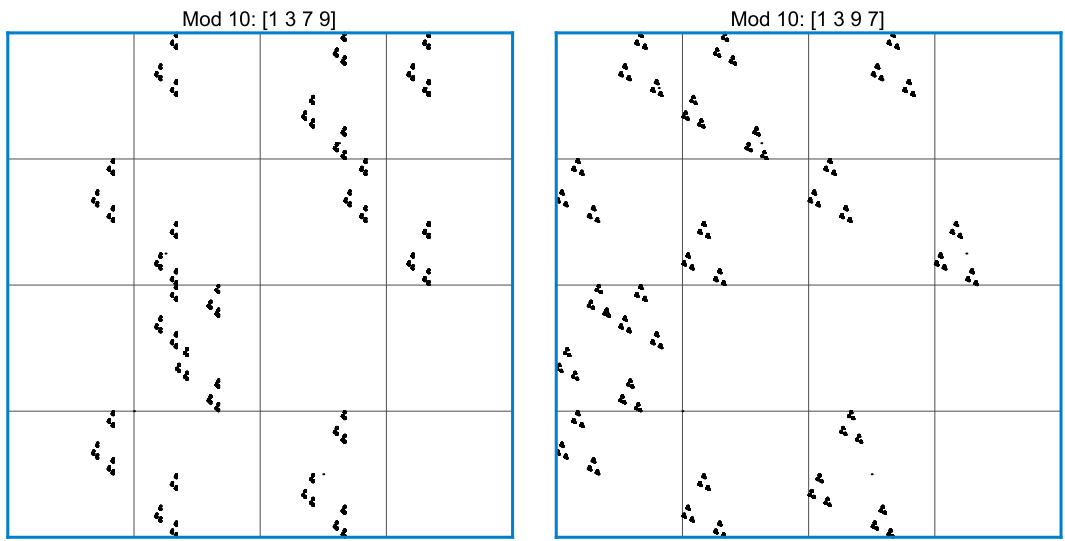}
\caption{IFS plot of residues modulo 10 of twin prime pairs (20,000 primes) beginning with $(5, 7)$.}
\label{fig10}
\end{figure}

\noindent 
What might this technique tell us about twin primes?  Fiqure \ref{fig10} shows the IFS plot for the residues modulo 10 of the sequence of twin prime pairs.  The first thing we notice is an abundance of forbidden addresses.  Table \ref{table7} shows why these empty addresses are expected.  For twin primes $(n, n+2)$, two things are always true: $n$ (mod 10) $\neq 3$ and $(n+2)$ (mod 10) $\neq 7$.  The reader can confirm in the left plot in Figure \ref{fig10}, that the addresses corresponding to the first five entries in Table \ref{table7} are empty. 

\vskip 6pt

Unlike the previous IFS plots, these plots are relatively sparse --- they approach their attractors very rapidly.  Indeed, after a few thousand points, there is little noticeable change.  Subdividing into longer addresses, we find that the number of empty addresses grows exponentially so that all of the plotted data is packed densely.  Nonetheless, we can see that addresses 14, 21, and 43 appear to be more populated than other length 2 addresses.  These correspond to the residue pairs $(9, 1), (1, 3)$, and $(7, 9)$, respectively.  Of the three, $(9,1)$ occurs most frequently because it is the only pair wherein both residues can appear either as a first or second element for twin primes.  Section 4 as well as \cite {LO-Sound} indicate that in the sequence of single primes, residue pair $(9, 1)$ occurs more frequently than $(1, 3)$ and $(7, 9)$.  For reference,

\vskip 6pt

\begin{center}
$\displaystyle \frac{\pi(7, 10^{8}; 10, (9,~1))}{\pi(7, 10^{8}; 10, (1,~3))} ~ \approx ~ \displaystyle \frac{\pi(7, 10^{8}; 10, (9,~1))}{\pi(7, 10^{8}; 10, (7,~9))} ~ \approx 1.08$. 
\end{center}

\begin{table}[htbp]
\footnotesize
\begin{center}
\begin{tabular}{| c | c | c | c |}
\hline Address & Residue Pair & Count  & Approx. \% \\ \hline
\hline
22 & $(3, 3)$ & 0 & 0\\ \hline
13 & $(7, 1)$ & 0 & 0\\ \hline
23 & $(7, 3)$ & 0 & 0\\ \hline
33 & $(7, 7)$ & 0 & 0\\ \hline
24 & $(9, 3)$ & 0 & 0\\ \hline
12 & $(3, 1)$ & 675,883 & 5.201\\ \hline
34 & $(9, 7)$ & 676,806 & 5.208\\ \hline
41 & $(1, 9)$ & 679,569 & 5.230\\ \hline
42 & $(3, 9)$ & 719,286 & 5.535\\ \hline
31 & $(1, 7)$ & 719,913 & 5.540\\ \hline
11 & $(1, 1)$ & 765,245 & 5.889\\ \hline
44 & $(9, 9)$ & 765,872 & 5.894\\ \hline
32 & $(3, 7)$ & 770,395 & 5.928\\ \hline
21 & $(1, 3)$ & 2,165,564 & 16.665\\ \hline
34 & $(7, 9)$ & 2,167,114 & 16.677\\ \hline
14 & $(9, 1)$ & 2,889,162 & 22.233\\ \hline
\end{tabular}
\caption{Frequencies for appearances of residue pairs modulo 10 for concatenation of individual members of twin primes up to $10^{8}$, sorted by count, starting with $(5, 7)$ (a total of 12,994,811 primes). See left plot in Figure \ref{fig10}.}
\label{table7}
\end{center}
\end{table}

\vspace*{-\baselineskip}

\begin{table}[htbp]
\footnotesize
\begin{center}
\begin{tabular}{| c | c | c | c |}
\hline Twin Pair Residues & $\pi_{\textmd{twin}}$\(11,10^{8}; 10, (a, b))$ & Percentage \\ \hline
\hline
$(1, 3)$ & 2,165,564 & 33.330\\ \hline
$(7, 9)$ & 2,167,114 & 33.354\\ \hline
$(9, 1)$ & 2,164,727 & 33.317\\ \hline
\end{tabular}
\caption{Frequencies for appearances of residue pairs modulo 10 for twin primes up to $10^{8}$.}
\label{table8}
\end{center}
\end{table}

\noindent
Looking at Table \ref{table8}, we see that, nonetheless, the three twin prime residue pairs occur with roughly the same frequency, and that, in fact, residue pair $(9, 1)$ seems to occur slightly \emph{less} frequently in the context of twin primes.

\vskip 6pt

Among the other entries in Table \ref{table7}, we can observe a bias at work.  For instance, the pair $(1, 3)$ is most often followed by $(7, 9)$ (address 32 in Table \ref{table7}). In comparison, the pair $(9, 1)$ follows $(1, 3)$ 93\% as often, and $(1, 3)$ follows itself only 88\% as often.  Indeed, looking at addresses 12, 34, and 41, we see that, just as with consecutive single primes, the residues of consecutive first members of twin primes do not like to repeat themselves.  It appears then that the prime biases observed in Section 4 do not favor one class of twin primes over the other two, but rather affect the order in which twin primes appear.

\vskip 6pt

What is surprising is to look at the values of $n$ for which $n\pm 1$ is prime, denoting them by $n_{\textmd{two}}$ (where ``two'' signifies ``2-tuple''). These values are necessarily of the form $n_{\textmd{two}} \equiv 0 \pmod{6}$.  We can therefore assign them to our four IFS transformations by taking $n_{\textmd{two}} \pmod{8}$ which has residue classes \{0, 2, 4, 6\}.

\vskip 6pt

\begin{figure}[htbp]
\centering
\epsfxsize = 204pt
\epsfbox{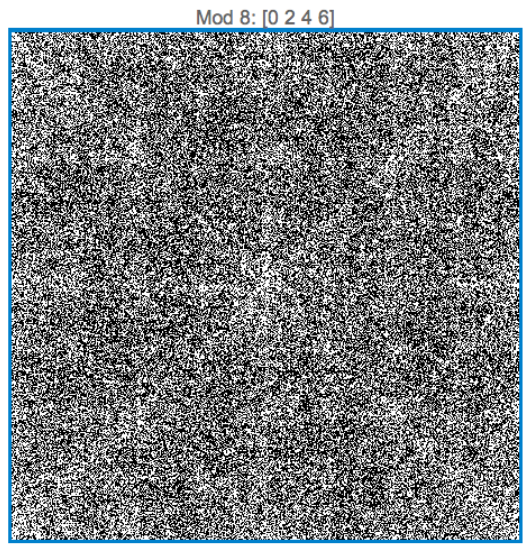}
\caption{IFS plot of residues modulo 8 of $n$ for which $n\pm 1$ is prime, $n$ up to $10^{6}$.}
\label{fig11}
\end{figure}

Although we are looking at the residues of \emph{composite} numbers, in Figure \ref{fig11} we see the ghostly suggestion of the IFS plots in Figure \ref{fig5}.  In Table \ref{table9}, we see that these composite numbers are indeed also biased and exhibit a similar repulsive behavior in that identical residues tend to avoid occurring in sequence.  In this range, the most common events occur roughly 1.13 times as often as the least common events.

\begin{table}[htbp]
\footnotesize
\begin{center}
\begin{tabular}{| c | c | c | c |}
\hline Address & Residue Pair & Count & Approx. \% \\ \hline
\hline
33 & $(4, 4)$ & 372,388 & 5.731\\ \hline
44 & $(6, 6)$ & 372,681 & 5.736\\ \hline
11 & $(0, 0)$ & 374,100 & 5.758\\ \hline
22 & $(2, 2)$ & 374,111 & 5.758\\ \hline
31 & $(0, 4)$ & 413,551 & 6.365\\ \hline
24 & $(6, 2)$ & 413,743 & 6.368\\ \hline
42 & $(2, 6)$ & 413,857 & 6.370\\ \hline
14 & $(4, 0)$ & 414,002 & 6.372\\ \hline
24 & $(4, 2)$ & 415,652 & 6.397\\ \hline
41 & $(0, 6)$ & 416,383 & 6.408\\ \hline
21 & $(2, 0)$ & 416,484 & 6.410\\ \hline
34 & $(6, 4)$ & 416,708 & 6.413\\ \hline
32 & $(2, 4)$ & 420,450 & 6.471\\ \hline
14 & $(6, 0)$ & 420,844 & 6.477\\ \hline
43 & $(4, 6)$ & 421,056 & 6.480\\ \hline
21 & $(0, 2)$ & 421,396 & 6.486\\ \hline
\end{tabular}
\caption{Length 2 addresses: frequencies for appearances of residues modulo 8 of $n$ for which $n\pm 1$ is prime, $n$ up to $10^{8}$.}
\label{table9}
\end{center}
\end{table}

\vspace*{-\baselineskip}

At first glance, it might appear that these composite numbers somehow encode information about the primes that immediately bracket them.  However, the same bias appears for $(n_{\textmd{two}}+k) \pmod{8}$ for \emph {any} 
value of $k$.  This follows directly from modular addition: although the frequency for each residue class can change depending on the value of $k$, the shape of the overall distribution remains the same.

\vskip 16pt

\section{Prime $k$-tuples}
\begin{figure}[htbp]
\centering
\epsfxsize = 420pt
\epsfbox{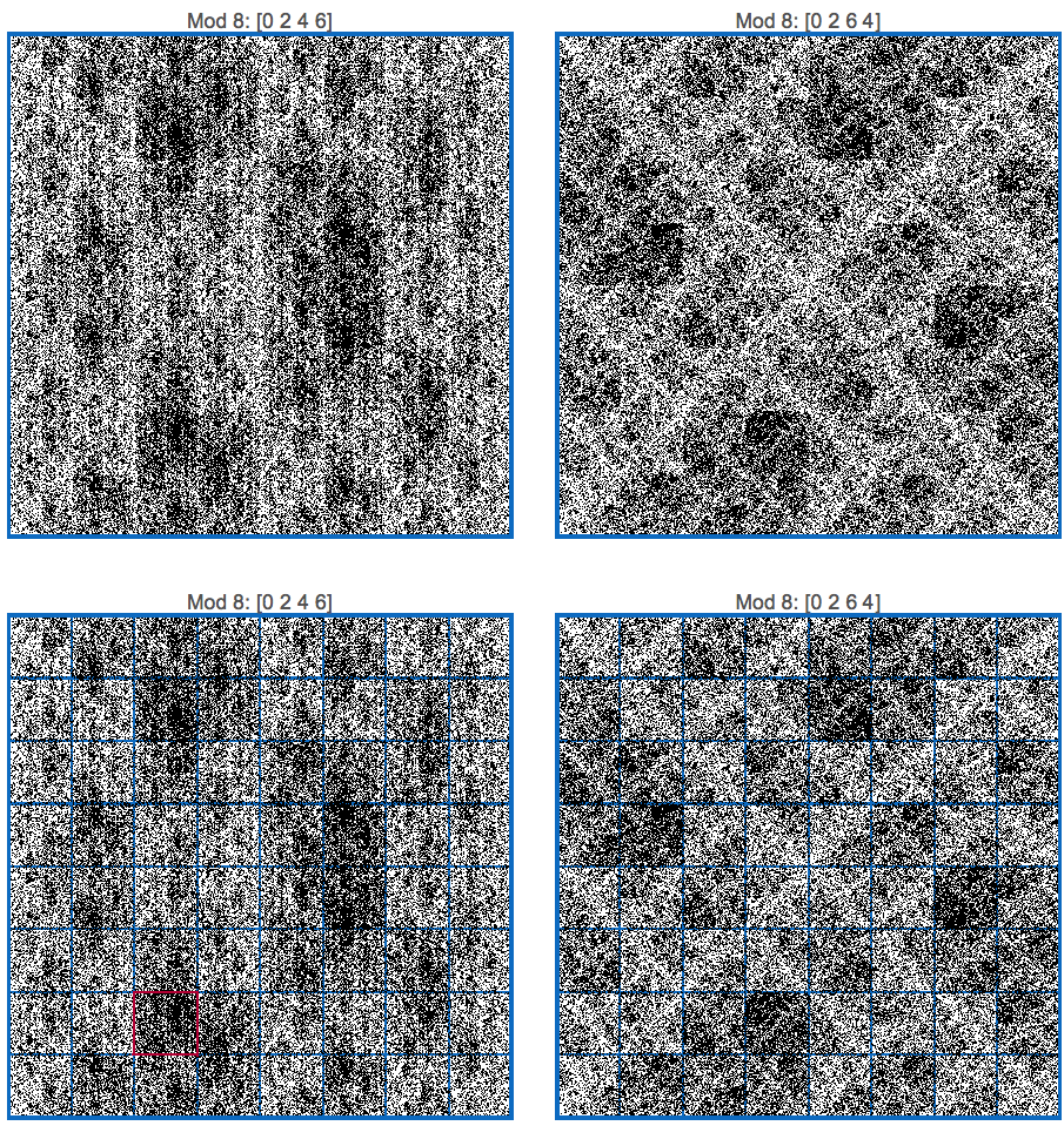}
\caption{IFS plot of residues modulo 8 of $n$ for which $n\pm 3$ is prime, $n$ up to $10^{6}$. The dark region highlighted in red in the lower left-hand plot is address 123.  Among length 3 addresses, the other darkest regions correspond to addresses 234, 341, and 412.  See Table \ref{table11} to zoom in to length 4 addresses.}
\label{fig12}
\end{figure}

\noindent 
These types of biases are not limited to twin primes, but rather appear to hold for other prime $k$-tuples, though not necessarily modulo 8.  For example, Figure \ref{fig12} and Table \ref{table10} show the results for $n_{\textmd{six}}$, associated with what are known as ``sexy primes."  In this range, the most common events occur roughly 1.47 times as often as the least common events.

\begin{table}[htbp]
\footnotesize
\begin{center}
\begin{tabular}{| c | c | c | c |}
\hline Address & Residue Pair & Count & Approx. \% \\ \hline
\hline
33 & $(4, 4)$ & 634,582 & 5.468\\ \hline
44 & $(6, 6)$ & 635,141 & 5.473\\ \hline
11 & $(0, 0)$ & 635,807 & 5.479\\ \hline
22 & $(2, 2)$ & 635,868 & 5.479\\ \hline
42 & $(2, 6)$ & 664,624 & 5.727\\ \hline
24 & $(6, 2)$ & 664,815 & 5.729\\ \hline
31 & $(0, 4)$ & 665,119 & 5.732\\ \hline
13 & $(4, 0)$ & 665,988 & 5.739\\ \hline
32 & $(2, 4)$ & 666,263 & 5.741\\ \hline
14 & $(6, 0)$ & 666,274 & 5.741\\ \hline
43 & $(4, 6)$ & 666,379 & 5.742\\ \hline
21 & $(0, 2)$ & 667,096 & 5.749\\ \hline
23 & $(4, 2)$ & 933,377 & 8.043\\ \hline
34 & $(6, 4)$ & 934,362 & 8.052\\ \hline
12 & $(2, 0)$ & 934,402 & 8.052\\ \hline
41 & $(0, 6)$ & 934,449 & 8.052\\ \hline
\end{tabular}
\caption{Length 2 addresses: frequencies for appearances of residues modulo 8 of $n$ for which $n\pm 3$ is prime, $n$ up to $10^{8}$.}
\label{table10}
\end{center}
\end{table}

We can tell at first glance that the distribution is far less heterogeneous than in the case of $n_{\textmd{two}}$. Furthermore, Table \ref{table11} points to a significant level of structure in the dynamics of $n_{\textmd{six}}$.  Indeed, they robustly resemble the findings in Section 5, in that the same residue class avoids appearing in sequence and the most frequent type of motion corresponds to stepping through the residue classes.  Here, the 4 most common sequences occur precisely in cyclical descending order through the residue classes.  In this range, these most frequent events occur roughly 5.47 times as often as the least common events.  Here again, we should emphasize that these are {\it composite} numbers exhibiting the repulsive characteristics that, up until now, have been associated with prime numbers.

\begin{table}[htbp]
\footnotesize
\begin{center}
\begin{tabular}{| c | c | c | c |}
\hline Address & Residue Quadruplet & Count & Approx. \% \\ \hline
\hline
2222 & $(2, 2, 2, 2)$ & 28,379 & 2.446\\ \hline
1111 & $(0, 0, 0, 0)$ & 28,517 & 2.457\\ \hline
3333 & $(4, 4, 4, 4)$ & 28,525 & 2.458\\ \hline
4444 & $(6, 6, 6, 6)$ & 28,815 & 2.483\\ \hline
2341 & $(0, 6, 4, 2)$ & 95,995 & 8.272\\ \hline
3412 & $(2, 0, 6, 4)$ & 96,074 & 8.279\\ \hline
4123 & $(4, 2, 0, 6)$ & 96,123 & 8.283\\ \hline
1234 & $(6, 4, 2, 0)$ & 96,189 & 8.289\\ \hline
\end{tabular}
\caption{Length 4 addresses: lowest and highest frequencies for appearances of residues modulo 8 of $n$ for which $n\pm 3$ is prime, $n$ up to $10^{8}$.}
\label{table11}
\end{center}
\end{table}

\FloatBarrier
\section{Conclusion}
\noindent 
Running the analysis from Section 4 with 1 million primes starting with the 4th prime, the millionth prime, or the 10 millionth prime shows similar results.  However, the plots are more heterogeneous and the histograms more evenly distributed for equal size samples of larger primes.  This suggests that the observed biases slowly even out as primes grow large (see Appendix A).  It is an open question as to whether, in the limit, these biases eventually give way to a uniform distribution \cite{LO-Sound}.

\vskip 6pt

Intriguing results can be obtained by applying appropriate sets of IFS transforms corresponding to the vertices of other equilateral shapes. Plots in 3 dimensions using $M_{15}$, $M_{16}$, $M_{20}$, or $M_{24}$ are compelling and can simultaneously embody different relationships depending on the plane from which they are viewed.  Preliminary evidence suggests that when examining $\pi(x_{0},x; q, (a, b))$, the self-avoidance bias for reduced residue classes holds for all $q>2$.

\vskip 6pt

At minimum, this work provides a compelling complement to the work of Lemke Oliver and Soundararajan.  The results with regard to dynamics, the distribution of prime $k$-tuples, and the self-avoiding characteristics of the residues of $k$-tuple-related composite numbers all appear to be new.  It is hoped that the use of IFS can provide further insight and avenues for exploration in the field of number theory.

\vskip 24pt

\noindent {\bf Acknowledgements.}  The author would like to thank Michael Frame and Michael Barnsley for their support and feedback.  Thanks also to Robert Schneider at Emory University for editorial and typesetting assistance.

\vskip 12pt

\noindent {\bf Postscript.} The author has come upon the colorful work of Chung-Ming Ko who performed a similar investigation in 2002, using a technique referred to as a ``two-dimensional histogram'' \cite{Ko}. Ko came to the same conclusions reached in Section 4 of the present paper (and the basic findings of Lemke-Oliver and Soundararajan), however the results were not interpreted from a dynamical standpoint and were not further developed.   

\vskip 20pt

\noindent \large{A}\small{PPENDICES}
\appendix

\section{Standard Deviation for Distribution\\of $\pi(x_{0},x_{0}+10^{6}; 10, (a, b))$ at Increasing Start Values $x_{0}$}
\vskip -10pt
\begin{table}[htbp]
\footnotesize
\begin{center}
\begin{tabular}{| c | c |}\hline
${x_{0}}$ & $\sigma$ (counts)\\ \hline \hline
7        & 15,640.6\\ \hline
$10^{6}$ & 14,366.5 \\ \hline
$10^{7}$ & 12,949.4\\ \hline
$10^{8}$ & 11,628.9\\ \hline
$10^{9}$ & 10,623.0\\ \hline
\hphantom{0}$10^{10}$& \hphantom{00}9,786.94\\ \hline
\hphantom{0}$10^{11}$& \hphantom{00}8,951.92\\ \hline
\hphantom{0}$10^{12}$& \hphantom{00}8,394.25\\ \hline
\end{tabular}
\caption{Standard deviation $\sigma$ is taken across counts for all $(a,b)$, where count=$\pi(x_{0},x_{0}+10^{6}; 10, (a, b))$.}
\label{append1}
\end{center}
\end{table}

\vspace*{-\baselineskip}

\section{Mathematica Code}
{\tiny
\begin{verbatim}
divider = 2;(* Set grid divisions for x and y *)
coPrime[k_] := Select[Range[k], CoprimeQ[#, k] == True &];
perm = Permutations[{1, 2, 3, 4}][[;; 2]];
ifsPlots := Module[{a, b, c, d},
  If[!MemberQ[{5, 8, 10, 12}, mod], {Print[Style["Please choose modulus 5, 8, 10, or 12", 14]], Abort[]}]
   Do[
    cp[i] = coPrime[mod][[perm[[i]]]];
    Thread[{a, b, c, d} = perm[[i]]];
    tranformOrder[i] = Switch[Mod[#, mod], coPrime[mod][[1]], a, coPrime[mod][[2]], b, coPrime[mod][[3]], c,
      coPrime[mod][[4]], d] & /@Prime[Range[start, start + span - 1]];
    f[j_, {x_, y_}] := 0.5*{x, y} + 0.5*Reverse[IntegerDigits[j - 1, 2, 2]]; 
    pt = {0.5, 0.5};
    ptlst = Table[pt = f[tranformOrder[i][[j]], pt], {j, Length[tranformOrder[i]]}];
    plts[i] = ListPlot[ptlst,
      Frame -> True,
      FrameStyle -> Directive[RGBColor[0, .5, .8], Thickness[.006]],
      FrameTicks -> None,
      GridLines -> {Table[k/divider, {k, 1, divider - 1}], 
      GridLinesStyle -> Directive[RGBColor[.3, .3, .3], Thickness[.002]],
      PlotRange -> {{0, 1}, {0, 1}},
      PlotStyle -> {PointSize -> ptSize, Black},
      AspectRatio -> 1,
      ImageSize -> imgSize, 
      PlotLabel -> Style["Mod " <> ToString@mod <> ": [" <> ToString@cp[i][[1]] <> " " <> ToString@cp[i][[2]] 
                  <> " " <> ToString@cp[i][[3]] <> " " <> ToString@cp[i][[4]] <> "]", 14]],
    {i, 2}];
  display[mod] = GraphicsRow[{plts[1], plts[2]}, Spacings -> 12]]

start = 4; span = 100000; imgSize = 360; ptSize = .001; mod = 10; ifsPlots

\end{verbatim}}

\end{document}